\documentclass{amsproc}
\usepackage{amsmath,amssymb,amsfonts}
\usepackage[mathscr]{eucal}
\usepackage{verbatim}

\usepackage{enumerate}

\numberwithin{equation}{section}

\theoremstyle{plain}
\newtheorem{theorem}[equation]{Theorem}

\newtheorem{prop}[equation]{Proposition}

\newtheorem{conj}[equation]{Conjecture}

\theoremstyle{definition}
\newtheorem{definition}[equation]{Definition}

\newtheorem*{thank}{Acknowledgments}

\input xy
\xyoption{all}

\newcommand{\Deltaop}{{\bf \Delta}^{op}}
\newcommand{\colim}{\text{colim}}

\newcommand{\nerve}{\text{nerve}}
\newcommand{\Hom}{\text{Hom}}
\newcommand{\SSets}{\mathcal{SS}ets}

\newcommand{\Sets}{\mathcal Sets}

\newcommand{\Secat}{\mathcal Se \mathcal Cat}
\newcommand{\Ho}{\text{Ho}}
\newcommand{\map}{\text{map}}
\newcommand{\Map}{\text{Map}}
\newcommand{\ob}{\text{ob}}

\newcommand{\css}{\mathcal{CSS}}

\newcommand{\Qcat}{\mathcal{QC}at}

\newcommand{\Thetanop}{\Theta_n^{op}}

\begin{document}

\title{Models for $(\infty,n)$-categories and the cobordism hypothesis}

\author[J.E. Bergner]{Julia E. Bergner}

\address{Department of Mathematics \\
University of California \\
Riverside, CA 92521}

\email{bergnerj@member.ams.org}

\date{\today}



\thanks{The author was partially supported by NSF grant DMS-0805951.}

\begin{abstract}
In this paper we introduce the models for $(\infty, n)$-categories which have been developed to date, as well as the comparisons between them that are known and conjectured.  We review the role of $(\infty, n)$-categories in the proof of the Cobordism Hypothesis.
\end{abstract}

\maketitle

\section{Introduction}

The role of higher categories is not new in the study of topological quantum field theories.  However, recent work of Lurie has introduced a homotopical approach to higher categories, that of $(\infty, n)$-categories, to the subject with his recent paper on the Cobordism Hypothesis.  The aim of this paper is to describe some of the known models for $(\infty, n)$-categories and the comparisons between them and to give a brief exposition of how they are used in Lurie's work.

For any positive integer $n$, an $n$-category consists of objects, 1-morphisms between objects, 2-morphisms between 1-morphisms, and so forth, up to $n$-morphisms between $(n-1)$-morphisms.  One can even have such higher morphisms for all $n$, leading to the idea of an $\infty$-category.  If associativity and identity properties are required to hold on the nose, then we have a strict $n$-category or strict $\infty$-category, and there is no problem with this definition.  However, in practice the examples that we find throughout mathematics are rarely this rigid.  More often we have associativity holding only up to isomorphism and satisfying some kinds of coherence laws, leading to the idea of a weak $n$-category or a weak $\infty$-category.  Many definitions have been proposed for such higher categories, but showing that they are equivalent to one another has proven to be an enormously difficult task.

Interestingly enough, the case of weak $\infty$-groupoids, where all morphisms at all levels are (weakly) invertible, can be handled more easily.  Given any topological space, one can think of it as an $\infty$-groupoid by regarding its points as objects, paths between the points as 1-morphisms, homotopies between the paths as 2-morphisms, homotopies between the homotopies as 3-morphisms, and continuing thus for all $k$-morphisms.  In fact, it is often taken as a definition that an $\infty$-groupoid \emph{is} a topological space.

If we take categories enriched in topological spaces, then we obtain a model for $(\infty, 1)$-categories, in which we have $k$-morphisms for all $k\geq 1$, but now they are only invertible for $k>1$; the points of the spaces now play the role of 1-morphisms.  Using the homotopy-theoretic equivalence between topological spaces and simplicial sets, it is common instead to consider categories enriched over simplicial sets, often called simplicial categories.  Topological or simplicial categories are good for many applications, but for others they are still too rigid, since composition of the mapping spaces is still required to have strict associativity and inverses.

There are different ways to weaken the definition of simplicial category so that composition is no longer strictly defined.  Segal categories, complete Segal spaces, and quasi-categories were all developed as alternatives to simplicial categories.  In sharp contrast to definitions of weak $n$-categories, these models are known to be equivalent to one another in a precise way.  More specifically, there is a model category corresponding to each model, and these model categories are Quillen equivalent to one another.  Further details on these models and the equivalences between them are given in the next section.

One might ask if we could move from $(\infty, 1)$-categories to more general $(\infty, n)$-categories using similar methods.  Not surprisingly, there are many more possible definitions of models for more general $n$.  In this paper, we describe several of these approaches as well as some of the known comparisons between them.  Many more of these relationships are still conjectural, although they are expected to be established in the near future.  Indeed, the treatment of different models in this paper is regrettably unbalanced, due to the fact that much of the work in this area is still being done.

In the world of topological quantum field theory, these kinds of higher-categorical structures have become important due to their important role in Lurie's recent proof of the Baez-Dolan Cobordism Hypothesis \cite{luriech}.  In the last section of this paper, we give a brief description of how to obtain a definition of a cobordism $(\infty, n)$-category.  We also give an introduction to the Cobordism Hypothesis as originally posed by Baez and Dolan and as proved by Lurie.

Throughout this paper we freely use the language of model categories and simplicial sets; readers unfamiliar with these methods are encouraged to look at \cite{ds} and \cite{gj} for further details.

\begin{thank}
Many thanks are due to the people who shared their knowledge of their work on this subject with me, including Clark Barwick, Jacob Lurie, Chris Schommer-Pries, and Claire Tomesch.  Helpful comments on the paper from the anonymous referee as well as editing suggestions from Arthur Greenspoon are also gratefully acknowledged.
\end{thank}

\section{Comparison of models for $(\infty, 1)$-categories}

In this section we review the various models for $(\infty, 1)$-categories, their model structures, and the Quillen equivalences between them.  A more extensive survey is given in \cite{survey}.

In some sense, the most basic place to start is with $(\infty, 0)$-categories, or $\infty$-groupoids.  Topological spaces or, equivalently, simplicial sets, are generally taken to be the definition of $\infty$-groupoids.  Certainly one can think of a topological space as an $\infty$-groupoid by regarding the points of the space as objects, the paths between points as 1-morphisms, homotopies between paths to be 2-morphisms, homotopies between homotopies as 3-morphisms, and so forth for all natural numbers.  While it can be argued that this explanation is one-directional, there is great difficulty in pinning down an actual precise definition of an $\infty$-groupoid, so it is common to take the above as a definition.

When moving one level higher to $(\infty, 1)$-categories, the most natural definition is to take topological categories, or categories enriched over topological spaces.  Now points of mapping spaces are 1-morphisms, paths are 2-morphisms, and so forth, so that now there is no reason for 1-morphisms to be invertible, but all higher morphisms are.  More commonly, we use \emph{simplicial categories}, or categories enriched over simplicial sets.

Given two objects $x$ and $y$ of a simplicial category $\mathcal C$, we denote the mapping space between them by $\Map(x,y)$.  For a simplicial category $\mathcal C$, its \emph{category of components} $\pi_0 \mathcal C$ is the ordinary category with objects the same as those of $\mathcal C$ and with objects given by
\[ \Hom_{\pi_0 \mathcal C}(x,y) = \pi_0 \Map_\mathcal C (x,y). \]
A simplicial functor $f \colon \mathcal C \rightarrow \mathcal D$ is a \emph{Dwyer-Kan equivalence} if
\begin{itemize}
\item for any objects $x, y$ of $\mathcal C$, the map
\[ \Map_{\mathcal C}(x,y) \rightarrow \Map_{\mathcal D}(fx,fy) \]
is a weak equivalence of simplicial sets, and

\item the induced functor on component categories
\[ \pi_0 f \colon \pi_0 \mathcal C \rightarrow \pi_0 \mathcal D \] is an equivalence of
categories.
\end{itemize}
A simplicial functor $f: \mathcal C \rightarrow \mathcal D$ is a \emph{fibration} if
\begin{itemize}
\item for any objects $x$ and $y$ in $\mathcal C$, the map
\[ \Hom_\mathcal C (x,y) \rightarrow \Hom_\mathcal D (fx,fy) \]
is a fibration of simplicial sets, and

\item for any object $x_1$ in $\mathcal C$, $y$ in $\mathcal
D$, and homotopy equivalence $e:fx_1 \rightarrow y$ in $\mathcal
D$, there is an object $x_2$ in $\mathcal C$ and homotopy
equivalence $d:x_1 \rightarrow x_2$ in $\mathcal C$ such that
$fd=e$.
\end{itemize}

\begin{theorem} \cite{simpcat}
There is a model structure $\mathcal{SC}$ on the category of small simplicial categories with weak equivalences the Dwyer-Kan equivalences and the fibrations as above.
\end{theorem}

However, there are other possible models for $(\infty, 1)$-categories.  One complaint that one might have about simplicial categories is that the requirement that the composition of mapping spaces be associative is too strong.  We are thus led to the definition of Segal categories

First, recall that a simplicial category $\mathcal C$ as we have defined it is a special case of a more general simplicial object in the category of small categories, where we assume that all face and degeneracy maps are the identity on the objects so that the objects form a discrete simplicial set.  We can take the simplicial nerve to obtain a simplicial space with 0-space discrete.  For any simplicial space $X$ recall that we can define Segal maps
\[ X_k \rightarrow \underbrace{X_1 \times_{X_0} \cdots \times_{X_0} X_1}_k. \]  In the case of a nerve of a simplicial category, these maps are isomorphisms of simplicial sets.  We obtain the notion of Segal category when we weaken this restriction on simplicial spaces.

\begin{definition} \cite{hs}
A simplicial space $X$ is a \emph{Segal precategory} if $X_0$ is discrete.  It is a \emph{Segal category} if, in addition, the Segal maps are weak equivalences of simplicial sets.
\end{definition}

In a Segal category $X$, we can consider the discrete space $X_0$ as the set of ``objects" and, given a pair of objects $(x,y)$, the ``mapping space" $\map_X(x,y)$ defined as the fiber over $(x,y)$ of the map $(d_1, d_0) \colon X_1 \rightarrow X_0 \times X_0$.  In this way, we can use much of the language of simplicial categories in the Segal category setting.  Furthermore, a Segal category $X$ has a corresponding homotopy category $\Ho(X)$ with the same objects as $X$ but with the sets of components of mapping spaces as the morphisms.

There is a functorial way of ``localizing" a Segal precategory to obtain a Segal category in such a way that the set at level zero is unchanged.  We denote this functor by $L$.  We then define a map $f \colon X \rightarrow Y$ of Segal precategories to be a \emph{Dwyer-Kan equivalence} if
\begin{enumerate}
\item for any $x$ and $y$ in $X_0$, the map $\map_{LX}(x,y) \rightarrow \map_{LY}(fx,fy)$ is a weak equivalence of spaces, and

\item the map $\Ho(LX) \rightarrow \Ho(LY)$ is an equivalence of categories
\end{enumerate}

\begin{theorem} \cite{thesis}, \cite{pel}
There are two model structures, $\Secat_c$ and $\Secat_f$, on the category of Segal precategories such that the fibrant objects are Segal categories and such that the weak equivalences are the Dwyer-Kan equivalences.
\end{theorem}

The need for two different model structures with the same weak equivalences arises in the comparison with other models; they are Quillen equivalent to one another via the identity functor.  In the model structure $\Secat_c$, the cofibrations are the monomorphisms and hence every object is cofibrant.  In $\Secat_f$, while it is not true that the fibrations are levelwise, the cofibrations are what one would expect them to be if they were; the discrepancy arises from technicalities in working with Segal precategories rather than with all simplicial spaces.

The following theorem can be regarded as a rigidification result, showing that weakening the condition on the Segal maps did not make much of a difference from the homotopy-theoretic point of view.

\begin{theorem} \cite{thesis}
The model categories $\mathcal{SC}$ and $\Secat_f$ are Quillen equivalent.
\end{theorem}

However, for many purposes the condition that the space at level zero be discrete is an awkward one.  Thus we come to our third model, that of complete Segal spaces.  These objects will again be simplicial spaces, and for technical reasons we require them to be fibrant on the Reedy model structure on simplicial spaces.

\begin{definition} \cite{rezk}
A Reedy fibrant simplicial space $W$ is a \emph{Segal space} if the Segal maps are weak equivalences of simplicial sets.
\end{definition}

Like simplicial categories, Segal spaces have objects (this time the set $W_{0,0}$) and mapping spaces between them.  Using the weak composition between mapping spaces, we can define homotopy equivalences and consider the subspace $W_h$ of such sitting inside of $W_1$.  It is not hard to see that the degeneracy map $W_0 \rightarrow W_1$ has image in $W_h$, since the image consists of ``identity maps" which are certainly homotopy equivalences.

\begin{definition} \cite{rezk}
A Segal space is \emph{complete} if the map $s_0 \colon W_0 \rightarrow W_h$ is a weak equivalence of simplicial sets.
\end{definition}

\begin{theorem} \cite{rezk}
There is a model structure $\css$ on the category of simplicial spaces in which the fibrant objects are the complete Segal spaces and the weak equivalences between fibrant objects are levelwise weak equivalences of simplicial sets.
\end{theorem}

\begin{theorem} \cite{thesis}
The model categories $\Secat_c$ and $\css$ are Quillen equivalent.
\end{theorem}

We approach the fourth model, that of quasi-categories, a bit differently.  If we begin with an ordinary category $\mathcal C$, its nerve is a simplicial set $\nerve(\mathcal C)$.  Again we can think of a Segal condition, here where the maps are isomorphisms of sets.  However, this description doesn't lend itself well to weakening, since we are dealing with sets rather than simplicial sets.  Alternatively, we can describe the ``composites" in the nerve of a category via what is commonly called a horn-filling condition.  Consider the inclusions $V[m,k] \rightarrow \Delta[m]$ for any $m \geq 1$ and $0 \leq k \leq m$.  A simplicial set $K$ is the nerve of a category if and only if any map $V[m,k] \rightarrow K$ extends uniquely to a map $\Delta[m] \rightarrow X$ for any $0<k<m$.  Because we don't include the cases where $k=0$ and $k=m$, this property is called the \emph{unique inner horn filling condition}.

In the special case where $\mathcal C$ is a groupoid, then $\nerve(\mathcal C)$ has the property that such a unique extension exists for $0 \leq k \leq m$, i.e., has the \emph{unique horn filling condition}.  However, in homotopy theory, it has long been common to consider simplicial sets that are a bit weaker than the nerves of groupoids.  If we have a simplicial set $K$ such that any map $V[m,k] \rightarrow K$ extends to a map $\Delta[m] \rightarrow K$, but this extension is no longer required to be unique, it is called a \emph{Kan complex}.  Such simplicial sets are significant in that they are the fibrant objects in the standard model structure on simplicial sets.  Therefore, they can be regarded as particular models for $\infty$-groupoids.

Here we return to the inner horn filling condition and call a simplicial set $K$ an \emph{inner Kan complex} or \emph{quasi-category} if it has the above non-unique extension property for $0<k<m$, a notion that was first defined by Boardman and Vogt \cite{bv}.  Then, just as a Kan complex is a homotopy version of a groupoid, a quasi-category is a homotopy version of a category, in fact a model for an $(\infty, 1)$-category.

\begin{theorem} \cite{dugspi}, \cite{joyalqcsc}, \cite{lurie}
There is a model structure $\Qcat$ on the category of simplicial sets such that the fibrant objects are the quasi-categories.
\end{theorem}

One can actually define mapping spaces, for example, in a quasi-category (\cite{dugspi} goes into particular detail on this point), and the weak equivalences between quasi-categories again have the same flavor as the Dwyer-Kan equivalences of simplicial sets.  Making this relationship more precise, there is a coherent nerve functor $\mathcal{SC} \rightarrow \Qcat$, first defined by Cordier and Porter \cite{cp}.  Proofs for using this functor to obtain a Quillen equivalence between $\mathcal{SC}$ and $\Qcat$ are given in Lurie's book and in unpublished work of Joyal; Dugger and Spivak have given a substantially shorter proof using different methods.

\begin{theorem} \cite{dugspi}, \cite{joyal}, \cite{lurie}
The model categories $\mathcal{SC}$ and $\Qcat$ are Quillen equivalent.
\end{theorem}

While the previous theorem was sufficient to prove that all four models are equivalent, Joyal and Tierney have established multiple direct Quillen equivalences between $\Qcat$ and the other models.

\begin{theorem} \cite{jt}
There are two different Quillen equivalences between $\css$ and $\Qcat$, and there are two analogous Quillen equivalences between $\Secat_c$ and $\Qcat$.
\end{theorem}

A fifth approach has long been of interest from the perspective of homotopy theory, namely that of viewing a simplicial category as a model for a homotopy theory, where the essential data of a ``homotopy theory" is a category with some specified class of weak equivalences.  This idea was made precise by Dwyer and Kan via their methods of simplicial localization \cite{dkfncxes}, \cite{dksimploc}.  Recent work of Barwick and Kan includes a model structure $\mathcal{CWE}$ on the category of small categories with weak equivalence, together with a Quillen equivalence between $\mathcal{CWE}$ and $\css$ \cite{bk}.

We should remark that these models are by no means the only ones which have been proposed; for example, $A_\infty$ categories are conjectured to be equivalent as well.

\section{Multisimplicial models: Segal $n$-categories and $n$-fold complete Segal spaces}

In this section we give definitions of the earliest defined models for $(\infty, n)$-categories.  The idea behind them is to iterate the simplicial structure, so an $(\infty, n)$-category is given by a functor $(\Deltaop)^n \rightarrow \SSets$, satisfying some properties.

\begin{definition}
An $n$-fold simplicial space is a functor $X \colon (\Deltaop)^n \rightarrow \SSets$.
\end{definition}

Notice that there are different ways to regard such an object as a functor.  One useful alternative is to think of an $n$-fold simplicial space as a functor
\[ X \colon \Deltaop \rightarrow \SSets^{(\Deltaop)^{n-1}} \]
where $\SSets^{(\Deltaop)^{n-1}}$ denotes the category of functors $(\Deltaop)^{n-1} \rightarrow \SSets$.  This perspective is useful in that it makes use of the idea that an $(\infty, n)$-category should somehow resemble a category enriched in $(\infty, n-1)$-categories.  In particular from this viewpoint we can consider the Reedy model structure.

The first definition for $(\infty, n)$-categories was that of Segal $n$-categories, first given by Hirschowitz and Simpson \cite{hs}.  It is given inductively, building from the definition of Segal category given in the previous section.  We denote by
\[ \SSets^{\Deltaop}_{disc} \] the category of Segal precategories, or functors $Y \colon \Deltaop \rightarrow \SSets$ such that $Y_0$ is discrete.  Then, define inductively the category $\SSets^{(\Deltaop)^n}_{disc}$ of functors $X \colon \Deltaop \rightarrow \SSets^{(\Deltaop)^{n-1}}_{disc}$ such that $X_0$ is discrete.  In particular, notice that discreteness conditions are built in at several levels.

\begin{definition}
An $n$-fold simplicial space $X \colon \Deltaop \rightarrow \SSets^{(\Deltaop)^{n-1}}_{disc}$ is a \emph{Segal} $n$-\emph{precategory} if $X_0$ is discrete.
It is a \emph{Segal} $n$-\emph{category} if, in addition, the Segal maps
\[ X_k \rightarrow \underbrace{X_1 \times_{X_0} \cdots \times_{X_0} X_1}_k \] are weak equivalences of Segal $(n-1)$-categories for $n \geq 2$.
\end{definition}

\begin{theorem} \cite{pel}
There is a model structure $n\Secat$ on the category of Segal $n$-precategories in which the fibrant objects are Segal $n$-categories.
\end{theorem}

To obtain a higher-order version of complete Segal spaces, we can work inductively, beginning with the definition of complete Segal spaces as given in the previous section.  Hence, in the following definitions we can assume that $n \geq 2$.  The definitions we give here are as stated by Lurie in \cite{luriech}; he gives a more general treatment of them in \cite{lurie2}.

\begin{definition}
A Reedy fibrant $n$-fold simplicial space is an $n$-\emph{fold Segal space} if
\begin{itemize}
\item each Segal map
\[ X_k \rightarrow \underbrace{X_1 \times_{X_0} \cdots \times_{X_0} X_1}_k \] is a weak equivalence of $(n-1)$-fold Segal spaces for each $k \geq 2$,
\item $X_k$ is an $(n-1)$-fold Segal space for each $k \geq 0$, and

\item the $(n-1)$-fold Segal space $X_0$ is essentially constant.
\end{itemize}
\end{definition}

Recall that an $(n-1)$-fold simplicial space $X$ is \emph{essentially constant} if there exists a weak equivalence $Y \rightarrow X$ where $Y$ is given by a constant diagram.

\begin{definition}
An $n$-fold Segal space $X$ is \emph{complete} if
\begin{itemize}
\item each $X_k$ is an $(n-1)$-fold complete Segal space, and

\item the simplicial space $X_{k, 0, \ldots, 0}$ is a complete Segal space for all $k \geq 0$.
\end{itemize}
\end{definition}

It is expected that there is a model category, which we denote $n\css$, on the category of $n$-fold simplicial spaces in which the fibrant objects are the $n$-fold complete Segal spaces.  Such a model structure seems to have been developed by Barwick but is not currently in the literature; it will be given precisely in \cite{inftyn2}.

\section{Models given by $\Theta_n$-diagrams}

As an alternative to $n$-fold complete Segal spaces, Rezk proposed a new model, that of $\Theta_n$-spaces \cite{rezktheta}.  The idea is to use a new diagram $\Theta_n$ rather than iterating simplicial diagrams.  We begin by defining the diagrams $\Theta_n$ inductively, using a more general construction on categories.

Given a category $\mathcal C$, define a category $\Theta \mathcal C$ with objects $[m](c_1, \ldots, c_m)$, where $[m]$ is an object of ${\bf \Delta}$ and $c_1, \ldots, c_m$ objects of $\mathcal C$.  A morphism
\[ [m](c_1, \ldots, c_m) \rightarrow [p](d_1, \ldots, d_p) \] is given by $(\delta, f_{ij})$, where $\delta \colon [m] \rightarrow [p]$ is a morphism in ${\bf \Delta}$ and $f_{ij} \colon c_i \rightarrow d_j$ is defined for every $1 \leq i \leq m$ and $1 \leq j \leq q$ where $\delta(i-1) < j \leq \delta (i)$ \cite[3.2]{rezktheta}.

Let $\Theta_0$ be the terminal category with one object and only the identity morphism.  Inductively define $\Theta_n = \Theta \Theta_{n-1}$.  Notice that $\Theta_1$ is just ${\bf \Delta}$.

One perspective on the objects of $\Theta_n$ is that they are ``basic" (strict) $n$-categories in the same way that objects of ${\bf \Delta}$ are ``basic" categories, in the sense that they encode the basic kinds of composites that can take place.  Therefore, if we take functors $\Thetanop \rightarrow \SSets$ and require conditions guaranteeing composition up to homotopy and some kind of completeness, we get models for $(\infty, n)$-categories.  In the case $n=1$, we obtain complete Segal spaces; for higher values of $n$, describing these conditions becomes more difficult but can be done in an inductive manner.

The underlying category for this model is $\SSets^{\Thetanop}$, the category of functors $\Thetanop \rightarrow \SSets$.  The model structure we want is obtained as a localization of the injective model structure on this category.

First, given an object $[m](c_1, \ldots, c_m)$ in $\Theta_n$, we obtain a corresponding ``simplex" $\Theta[m](c_1, \ldots, c_m)$ in the category $\Sets^{\Thetanop}$; making it constant in the additional simplicial direction gives an object of $\SSets^{\Thetanop}$.  This object should be regarded as the analogue of the $m$-simplex $\Delta[m]$ arising from the object $[m]$ of ${\bf \Delta}$.

Given $m \geq 2$ and $c_1, \ldots, c_m$ objects of $\Theta_{n-1}$, define the object
\[ G[m](c_1, \ldots, c_m)= \colim (\Theta[1](c_1) \leftarrow \Theta[0] \rightarrow \cdots \leftarrow \Theta[0] \rightarrow \Theta[1](c_m)). \]  There is an inclusion map
\[ se^{(c_1, \ldots, c_m)} \colon G[m](c_1, \ldots, c_n) \rightarrow \Theta[n](c_1, \ldots, c_m). \]  We define the set
\[ Se_{\Theta_n} = \{ se^{(c_1, \ldots, c_m)} \mid m \geq 2, c_1, \ldots c_m \in \ob(\Theta_{n-1})\}. \]  Localizing with respect to this set of maps gives composition up to homotopy, but only on the level of $n$-morphisms.  We need to localize additionally in such a way that lower-level morphisms also have this property, and we can do so inductively.

In \cite[4.4]{rezktheta}, Rezk defines an intertwining functor
\[ V \colon \Theta(\SSets^{\Theta_{n-1}^{op}}_c) \rightarrow \SSets^{\Thetanop}_c \] by
\[ V[m](A_1, \ldots, A_m)([q](c_1, \ldots, c_q))= \coprod_{\delta \in \Hom_{\bf \Delta}([q], [m])} \prod_{i=1}^q \prod_{j=\delta(i-1)+1}^{\delta(i)} A_j(c_i) \] where the $A_j$ are objects of $\SSets^{\Theta_{n-1}^{op}}$ and the $c_i$ are objects of $\Theta_n$.  This functor can be used to ``upgrade" sets of maps in $\SSets^{\Theta_{n-1}^{op}}$ to sets of maps in $\SSets^{\Thetanop}$.  Given a map $f \colon A \rightarrow B$ in $\SSets^{\Theta_{n-1}^{op}}$, we obtain a map $V[1](f) \colon V[1](A) \rightarrow V[1](B)$.

Let $\mathcal S_1=Se_{\bf \Delta}$, and for $n \geq 2$, inductively define $\mathcal S_n=Se_{\Theta_n} \cup V[1](\mathcal S_{n-1})$.   Localizing the model structure $\SSets^{\Thetanop}_c$ with respect to $\mathcal S_n$ results in a cartesian model category whose fibrant objects are higher-order analogues of Segal spaces.

However, we need to incorporate higher-order completeness conditions as well.  To define the maps with respect to which we need to localize, we make use of an adjoint relationship with simplicial spaces as described by Rezk in \cite[4.1]{rezktheta}.  First, define the functor $T \colon {\bf \Delta} \rightarrow \SSets^{\Thetanop}$ by
\[ T[q]([m](c_1, \ldots, c_m)) = \Hom_{\bf \Delta}([m], [q]). \]  We use this functor $T$ to define the functor $T^* \colon \SSets^{\Thetanop} \rightarrow \SSets^{\Deltaop}$ defined by
\[ T^*(X)[m]= \Map_{\SSets^{\Thetanop}}(T[m], X), \] which has a left adjoint $T_\#$.  This adjoint pair is in fact a Quillen pair with respect to the injective model structures.

Now, define $Cpt_{\bf \Delta}= \{E \rightarrow \Delta[0]\}$ and, for $n \geq 2$,
\[ Cpt_{\Theta_n}=\{T_\# E \rightarrow T_\# \Delta[0]\}. \]  Let $\mathcal T_1=Se_{\Theta_1} \cup Cpt_{\Theta_1}$ and, for $n \geq 2$,
\[ \mathcal T_n=Se_{\Theta_n} \cup Cpt_{\Theta_n} \cup V[1](\mathcal T_{n-1}). \]

\begin{theorem} \cite[8.5]{rezktheta}
Localizing $\SSets^{\Thetanop}_c$ with respect to the set $\mathcal T_n$ gives a cartesian model category, which we denote by $\Theta_nSp$.  Its fibrant objects are higher-order analogues of complete Segal spaces, and therefore it is a model for $(\infty, n)$-categories.
\end{theorem}

In his work on the categories $\Theta_n$, Joyal suggested that there should be models for $(\infty, n)$-categories given by functors $\Theta_n \rightarrow \Sets$ satisfying higher-order inner horn-filling conditions, but he was unable to find the right way to describe these conditions.  More recently, Barwick has been able to use the relationship between complete Segal spaces and quasi-categories to formulate a higher-dimensional version.

\begin{definition} \cite{bar}
A \emph{quasi}-$n$-\emph{category} is a functor $\Theta_n^{op} \rightarrow \Sets$ satisfying appropriate inner horn-filling conditions.
\end{definition}

\begin{theorem} \cite{bar}
There is a cartesian model structure $n\Qcat$ on the category $\Sets^{\Thetanop}$ in which the fibrant objects are quasi-$n$-categories.
\end{theorem}

\section{$(\infty, n)$-categories as enriched categories: strict and weak versions}

Intuitively, one would like to think of $(\infty, n)$-categories as categories enriched over $(\infty, n-1)$-categories.  In practice, this approach can be problematic.  In particular, if we want our models for $(\infty, n)$-categories to be objects in a model category, then at the very least we need our model structure on $(\infty, n-1)$-categories to be cartesian.  Since the model structure $\mathcal{SC}$ for simplicial categories is not cartesian, in that the product is not compatible with the model structure, we cannot continue the induction using that model.  However, we can enrich over other models.

The model category $\Theta_nSp$ is cartesian, and therefore we can consider categories enriched in $\Theta_{n-1}Sp$ as another model for $(\infty, n)$-categories.  In doing so, we have a way to realize the intuitive idea that $(\infty, n)$-categories are categories enriched in $(\infty, n-1)$-categories.

Let $\Map_\mathcal C(x,y)$ denote the mapping object in $\Theta_{n-1}Sp$ between objects $x$ and $y$ of a category $\mathcal C$ enriched in $\Theta_{n-1}Sp$.

\begin{definition} \cite{inftyn1}
Let $\mathcal C$ and $\mathcal D$ be categories enriched in $\Theta_{n-1}Sp$.  An enriched functor $f \colon \mathcal C \rightarrow \mathcal D$ is a \emph{weak equivalence} if
\begin{enumerate}
\item $\Map_{\mathcal C}(x,y) \rightarrow \Map_{\mathcal D}(fx,fy)$ is a weak equivalence in $\Theta_{n-1}Sp$ for any objects $x,y$, and

\item $\pi_0 \mathcal C \rightarrow \pi_0 \mathcal D$ is an equivalence of categories, where $\pi_0 \mathcal C$ has the same objects as $\mathcal C$ and $\Hom_{\pi_0 \mathcal C} (x,y) = \Hom_{\Ho(\Theta_{n-1}Sp)}(1, \Map_\mathcal C(x,y))$.
\end{enumerate}
\end{definition}

\begin{theorem} \cite{inftyn1}
There is a model structure $\Theta_{n-1}Sp-Cat$ on the category of small categories enriched in $\Theta_{n-1}Sp$ with weak equivalences defined as above.
\end{theorem}

Just as we had a simplicial nerve functor taking a simplicial category to a simplicial diagram of simplicial sets, we have a nerve functor taking a category enriched in $\Theta_{n-1}Sp$ to a simplicial diagram of objects in $\Theta_{n-1}Sp$.  If we call the resulting simplicial object $X$, we can observe that the strict Segal condition holds in this setting, in that the maps
\[ X_k \rightarrow \underbrace{X_1 \times_{X_0} \cdots \times_{X_0} X_1}_k \] are isomorphisms of objects in $\Theta_{n-1}Sp$.  Weakening this condition as before leads to the following definition.

\begin{definition} \cite{inftyn1}
A \emph{Segal precategory object in} $\Theta_{n-1}Sp$ is a functor $X \colon \Deltaop \rightarrow \Theta_{n-1}Sp$ such that the functor $X_0 \colon \Theta_{n-1}^{op} \rightarrow \SSets$ is discrete.  It is a \emph{Segal category object} if the Segal maps are weak equivalences in $\Theta_{n-1}Sp$.
\end{definition}

We can again define a functor $L$ taking a Segal precategory object to a Segal category object, and define mapping objects (now $\Theta_{n-1}$-spaces rather than simplicial sets) and a homotopy category for a Segal category object such as we did for a Segal category.  A functor of Segal precategory objects $f \colon X \rightarrow Y$ is a \emph{Dwyer-Kan equivalence} if
\begin{enumerate}
\item for any $x$ and $y$ in $X_0$, the map $\map_{LX}(x,y) \rightarrow \map_{LY}(fx,fy)$ is a weak equivalence in $\Theta_{n-1}Sp$, and

\item the map $\Ho(LX) \rightarrow \Ho(LY)$ is an equivalence of categories.
\end{enumerate}

\begin{theorem} \cite{inftyn1}
There are two model structures on the category of Segal precategory objects, $\mathcal Se(\Theta_{n-1}Sp)^{\Deltaop}_{disc,f}$ and $\mathcal Se(\Theta_{n-1}Sp)^{\Deltaop}_{disc, c}$, in which the weak equivalences are the Dwyer-Kan equivalences.
\end{theorem}

The notion of more general Segal category objects has also been developed by Simpson \cite{simpsonbook}.

\section{Comparisons between different models}

With so many models for $(\infty, n)$-categories, we need to establish that they are all equivalent to one another.  At this time, we have partial results in this area and many more conjectures.

We begin with known results, which can be summarized by the following diagram:
\[ (\Theta_{n-1}Sp)-Cat \leftrightarrows \mathcal Se(\Theta_{n-1}Sp)^{\Deltaop}_{disc, f} \rightleftarrows \mathcal Se(\Theta_{n-1}Sp)^{\Deltaop}_{disc, c}. \]

The leftmost Quillen equivalence is given by the following theorem, which is a generalization of the Quillen equivalence between $\mathcal{SC}$ and $\Secat_f$.

\begin{theorem} \cite{inftyn1}
The enriched nerve functor is the right adjoint for a Quillen equivalence between $(\Theta_{n-1}Sp)-Cat$ and $\mathcal Se(\Theta_{n-1}Sp)^{\Deltaop}_{disc, f}$.
\end{theorem}

Moving to the right, the next Quillen equivalence is the easiest.

\begin{prop} \cite{inftyn1}
The identity functor gives a Quillen equivalence between $\mathcal Se(\Theta_{n-1}Sp)^{\Deltaop}_{disc, f}$ and $\mathcal Se(\Theta_{n-1}Sp)^{\Deltaop}_{disc, c}$.
\end{prop}

To continue the comparison, we make use of a definition of \emph{complete Segal space objects in} $(\Theta_{n-1}Sp)^{\Deltaop}$.  There is a model structure $L_{CS}(\Theta_{n-1}Sp)^{\Deltaop}_c$ on the category of functors $\Deltaop \rightarrow \Theta_{n-1}Sp$ in which the fibrant objects satisfy Segal and completeness conditions.  The Segal condition is as given for Segal category objects above, but the completeness condition is subtle, and the comparison with Segal category objects is still work in progress \cite{inftyn2}.

\begin{conj}
There are Quillen equivalences
\[ \mathcal Se(\Theta_{n-1}Sp)^{\Deltaop}_{disc, c} \rightleftarrows L_{CS}(\Theta_{n-1}Sp)^{\Deltaop} \rightleftarrows \Theta_nSp. \]
\end{conj}

The second equivalence in this chain should in fact be the first in a chain of Quillen equivalences between $\Theta_nSp$ and $n\css$, induced by the chain of functors
\[ {\bf \Delta}^n = {\bf \Delta}^{n-1} \times {\bf \Delta} \rightarrow {\bf \Delta}^{n-1} \times \Theta_2 \rightarrow \cdots \rightarrow {\bf \Delta} \times \Theta_{n-1} \rightarrow \Theta_n. \]

\begin{conj} \label{chain}
There is a chain of Quillen equivalences
\[ n\css \rightleftarrows L_{CS}(\Theta_2Sp)^{(\Deltaop)^{n-2}} \rightleftarrows \cdots \rightleftarrows L_{CS}(\Theta_{n-1}Sp)^{\Deltaop} \rightleftarrows \Theta_nSp. \]
\end{conj}

There are other models, and proposed comparisons between them, that are currently being developed.  Since the model structure $n\Qcat$ is also cartesian, we can define categories enriched in $(n-1)\Qcat$.  With an appropriate model structure on the category of such, we conjecture that it is Quillen equivalent to $(\Theta_{n-1}Sp)-Cat$, using the Quillen equivalence between $\Theta_nSp$ and $n\Qcat$.

It is expected that there is yet another method of connecting $\mathcal Se(\Theta_{n-1}Sp)^{\Deltaop}_{disc, c}$ to $\Theta_nSp$ with Quillen equivalences, where the intermediate model category has objects functors $\Thetanop \rightarrow \SSets$ satisfying discreteness conditions.  In other words, such a model structure would be a Segal category version of $\Theta_nSp$.  This model structure should be connected to $n\Secat$ via a chain of Quillen equivalences analogous to those connecting $\Theta_nSp$ and $n\css$.  The comparison between $n\Secat$ and other models is also being investigated by Tomesch \cite{tom}.

\section{$(\infty, n)$-categories and the cobordism hypothesis}

In this section, we give a brief account of the Cobordism Hypothesis, which was originally posed in the context of weak $n$-categories by Baez and Dolan in \cite{bd} and proved in the setting of $(\infty, n)$-categories by Lurie \cite{luriech}.  This result not only gives a purely algebraic description of higher categories defined in terms of cobordisms of manifolds, but is usually interpreted as a statement about topological quantum field theories.

We begin with Atiyah's definition of topological quantum field theory \cite{atiyah}.

\begin{definition}
For $n \geq 1$, denote by $Cob(n)$ the category with objects closed framed $(n-1)$-dimensional manifolds and morphisms the diffeomorphism classes of framed cobordisms between them.
\end{definition}

Notice that composition is defined by gluing together cobordisms; since we are taking diffeomorphism classes on the level of morphisms there is no difficulty defining composition.  This category is also equipped with a symmetric monoidal structure given by disjoint union of $(n-1)$-dimensional manifolds.

\begin{definition}
For a field $\Bbbk$, a \emph{topological quantum field theory} of dimension $n$ is a symmetric monoidal functor $Z \colon Cob(n) \rightarrow Vect(\Bbbk)$, where $Vect(\Bbbk)$ is the category of vector spaces equipped with the usual tensor product.
\end{definition}

An important observation about topological quantum field theories is that the structure of $Cob(n)$ affects which vector spaces can be in the image of one.  Since a framed manifold has a dual, the vector space to which it is assigned must also have a well-behaved dual; in particular it must be finite-dimensional.

In low dimensions, topological quantum field theories are well-understood, since the manifolds appearing in the category $Cob(n)$ are easily described.  In particular, one can cut them into simpler pieces in a nice way and understand the entire functor just by knowing what happens to these pieces.  However, when $n \geq 3$, these pieces may not just be manifolds with boundary but instead manifolds with corners.  We need a higher-categorical structure to encode this larger range of dimensions of manifolds.

\begin{definition}
For $n \geq 1$, define a weak $n$-category $Cob_n(n)$ to have objects 0-dimensional closed framed 0-dimensional manifolds, 1-morphisms framed cobordisms between them, 2-morphisms framed cobordisms between the cobordisms, up to $n$-morphisms which are diffeomorphism classes of framed cobordisms of dimension $n$.
\end{definition}

Notice that, unlike for $Cob(n)$, we cannot take diffeomorphism classes except at the top dimension, since cobordisms must be defined between actual manifolds, not diffeomorphism classes of them.  Therefore, composition of morphisms is no longer strictly defined, and hence we have a weak, rather than strict, $n$-category.

To define a generalized topological field theories using $Cob_n(n)$ rather than $Cob(n)$, we need a higher-categorical version of $Vect(\Bbbk)$.  While there are several proposed definitions for such a weak 2-category, it is unknown how to continue to still higher dimensions.  Fortunately, we do not really need a particular structure and therefore can replace $Vect(\Bbbk)$ with an arbitrary symmetric monoidal weak $n$-category.

\begin{definition}
Let $\mathcal C$ be a symmetric monoidal weak $n$-category.  An \emph{extended} $\mathcal C$-\emph{valued topological field theory of dimension} $n$ is a symmetric monoidal functor $Z \colon Cob_n(n) \rightarrow \mathcal C$.
\end{definition}

We are almost ready to state Baez and Dolan's original conjecture, but we first need to comment on duality in higher categories.  As mentioned above, the image of a topological quantum field theory had to consist of finite-dimensional vector spaces, since they have well-behaved duals.  In the extended world, duality becomes more complicated, in that we need to consider objects as well as morphisms of all levels.  Objects having the appropriate properties are called \emph{fully dualizable objects}.  We do not give a precise definition of fully dualizable here, but refer the reader to Lurie's paper \cite[2.3]{luriech}.

\begin{theorem} (Baez-Dolan Cobordism Hypothesis)
Let $\mathcal C$ be a symmetric monoidal weak $n$-category and $Z$ a topological quantum field theory.  Then the evaluation functor $Z \mapsto Z(\ast)$ determines a bijection between isomorphism classes of framed extended $\mathcal C$-valued topological quantum field theories and isomorphism classes of fully dualizable objects of $\mathcal C$.
\end{theorem}

This version of the Cobordism Hypothesis has been proved in the case where $n=2$ by Schommer-Pries \cite{sp}, but it seems difficult to prove for higher values of $n$.  Not least is the difficulty of knowing how to handle weak $n$-categories, and how strict or weak the symmetric monoidal structure should be.  However, Lurie's main insight, allowing a sidestep to this problem, was to prove an alternative version using $(\infty, n)$-categories.  The original version follows via truncation methods.

Here, we focus on the definition of an $(\infty, n)$-category $Bord_n^{fr}$ which generalizes $Cob_n(n)$.  Following Lurie \cite{luriech}, we use the $n$-fold complete Segal space model.  An alternative definition in the setting of $\Theta_n$-spaces is being developed by Rozenblyum and Schommer-Pries \cite{rsp}.

\begin{definition} \cite[2.2]{luriech}
Let $V$ be a finite-dimensional vector space.  Define an $n$-fold simplicial space $PBord_n^{fr, V}$ by
\[ (PBord_n^{fr, V})_{k_1, \ldots, k_n}= \{(M, \{t_0^1 \leq \cdots, \leq t^1_{k_1}\}, \ldots, \{t_0^1 \leq \cdots, t_{k_n}^n\}) \} \] where:
\begin{itemize}
\item $M$ is a closed $n$-dimensional framed submanifold of $V \times \mathbb R^n$ such that the projection $M \rightarrow \mathbb R^n$ is proper,
\item for every $S \subseteq \{1, \ldots, n\}$ and every $\{0 \leq j_i \leq k_i \}_{i \in S}$, the composite map
\[ M \rightarrow \mathbb R^n \rightarrow \mathbb R^S \] does not have $(t_{j_i})_{i \in S}$ as a critical value, and

\item the map $M \rightarrow \mathbb R^{\{i+1, \ldots, n\}}$ is submersive at every $x \in M$ whose image in $\mathbb R^{\{i\}}$ belongs to the set $\{t_{i_0}, \ldots, t_{i_k}\}$.
\end{itemize}
As $V$ ranges over all finite-dimensional subspaces of $\mathbb R^\infty$, let
\[ PBord_n^{fr} = \colim_V PBord^{fr,V}_n. \]
\end{definition}

\begin{prop}
The $n$-fold simplicial space $PBord_n^{fr}$ is an $n$-fold Segal space which is not necessarily complete.
\end{prop}

To define $Bord_n^{fr}$, we use the fact that any $n$-fold simplicial space has a completion, or weakly equivalent $n$-fold complete Segal space.

\begin{definition}
Define $Bord_n^{fr}$ to be the completion of $PBord_n^{fr}$.
\end{definition}

With this definition in place, we can now state Lurie's version of the Cobordism Hypothesis.

\begin{theorem} \cite[1.4.9]{luriech}
Let $\mathcal C$ be a symmetric monoidal $(\infty, n)$-category.  The evaluation functor $Z \mapsto Z(\ast)$ determines a bijection between isomorphism classes of symmetric monoidal functors $Bord_n^{fr} \rightarrow \mathcal C$ and isomorphism classes of fully dualizable objects of $\mathcal C$.
\end{theorem}

We conclude with some comments about Lurie's method of proof, in particular his use of $(\infty, n)$-categories.  Using intuitive definitions, it would seem that $(\infty, n)$-categories would be more complicated than weak $n$-categories and therefore more difficult to use in practice.  However, the homotopical nature of $(\infty, n)$-categories makes the opposite true.  While in each case composition is defined weakly, for $(\infty, n)$-categories the resulting coherence data is nicely packaged into the definition.  This way of dealing with the inherent complications of weak $n$-categories proved to be an effective technique in establishing this result.

\end{document}